\documentclass[12pt]{article}  

\usepackage{hyperref}
\usepackage{graphicx}
\usepackage{amsfonts}
\usepackage{amsmath}
\usepackage{tagging}

\usepackage{enumerate}

\newtheorem{theorem}{Theorem}

\newtheorem{proposition}[theorem]{Proposition}
\newtheorem{lemma}[theorem]{Lemma}
\newtheorem{corollary}[theorem]{Corollary}
\newtheorem{remark}[theorem]{Remark}
\newtheorem{example}[theorem]{Example}

\newcommand{\proof}{ \noindent{\it Proof:\ \ }}
\def\qed{\ifhmode\unskip\nobreak\fi\ifmmode\ifinner\else\hskip5 pt
\fi\fi\hbox{\hskip5 pt \vrule width4 pt height6 pt depth1.5 pt
\hskip 1pt }}
\newcommand{\po}{{\hspace*{-1ex}}{\bf .  }}

\let\oldmarginpar\marginpar
\renewcommand\marginpar[1]{${}^\clubsuit$\oldmarginpar[\raggedleft\scriptsize\sf #1]{\raggedright\scriptsize\sf #1}}
\newcommand{\cref}[1]{Corollary~\ref{#1}}

\newcommand{\lref}[1]{Lemma~\ref{#1}}
\newcommand{\pref}[1]{Proposition~\ref{#1}}
\newcommand{\rref}[1]{Remark~\ref{#1}}
\newcommand{\sref}[1]{Section~\ref{#1}}
\newcommand{\tref}[1]{Theorem~\ref{#1}}

\def\Bbb#1{\mathbb#1}
\newcommand{\R}{\Bbb{R}}
\newcommand{\Z}{\Bbb{Z}}
\newcommand{\Q}{\Bbb{Q}}
\newcommand{\E}{\Bbb{E}}

\newcommand{\RP}{\Bbb{R}\Bbb{P}}
\newcommand{\Sp}{\Bbb{S}}
\newcommand{\DS}{\Bbb{S}_{-1}}
\newcommand{\DSc}{\Bbb{S}_c}

\newcommand{\Hi}{\Bbb{H}}

\newcommand{\RR}{\Bbb{R}_+^3}

\newcommand{\ra}{\rangle}
\newcommand{\la}{\langle}

\newcommand{\vs}{\vspace{.25in}}
\newcommand{\vsss}{\vspace{.1in}}

\newcommand{\spa}{\mbox{span}}
\newcommand{\im}{{\rm Im \,}}

\newcommand{\sinc}{S}
\newcommand{\cosc}{C}

\newcommand{\tr}{{\rm tr\,}}
\renewcommand{\d}{\mathrm{d}}
\newcommand{\SO}{\operatorname{SO}}

\usetag{pics}
\newcommand*\checkusemytag
{
    \begingroup
    \def\tempa{nopics}
    \expandafter
    \endgroup
    \ifx\tempa\usemytag
        \droptag{pics}
    \else
    \fi
}
\checkusemytag

\begin{document}

\title{Curvature homogeneous hypersurfaces\\in space forms}  
\author{Robert Bryant, Luis Florit and Wolfgang Ziller
\thanks{The first author was supported by a grant from the Simons
Foundation (347349, Bryant). The third author was supported by a ROG
grant from the University of Pennsylvania and would like to thank IMPA
for its hospitality.}
}
\date{}
\maketitle

\section{Introduction} 

Since the early days of differential geometry we know that a metric on a
manifold defines a natural notion of curvature, collected in what today
is called its curvature tensor. In particular, isometries preserve this
curvature, and thus the curvature of a homogeneous space should be
somehow constant along the manifold. This logically leads to the question
of whether a Riemannian manifold whose curvature tensor is “the same” at
every point is in fact homogeneous. More generally, a natural problem in
Riemannian geometry is to what extent, and even in which sense, “the
curvature determines the metric”. This question is more subtle than it
seems at first glance, having several answers and open aspects.

More precisely, a Riemannian manifold $M^n$ is called {\it curvature
homogeneous} if, for any pair of points $p,q\in M$, there exists a linear
isometry $J_{pq}:T_pM\to T_qM$ that preserves the curvature tensor $R$,
i.e., $J_{pq}^* R_q=R_p$. In dimension $n=3$ this is equivalent to the
condition that the eigenvalues of the Ricci tensor are constant.
I. Singer showed in \cite{Si} that, if sufficiently many covariant
derivatives of $R$ match, then the metric is in fact homogeneous. He then
asked whether it is sufficient to only require this condition for $R$.
Since then, this has been an active research area and the question turns
out to be surprisingly subtle.

The first complete non-homogeneous examples were given by K. Sekigawa and
H. Takagi in \cite{Se,Ta}, where they showed that such examples exist in
any odd dimension and depend on several functions of one variable. On the
other hand, there are also several known obstructions. In \cite{TV} it
was shown that, if the curvature tensor $R$ is that of an irreducible
symmetric space, then the metric itself must be locally symmetric. In the
above examples, $R$ was that of $\Hi^2\times\R^{n-2}$ and these were in
fact classified in \cite{Bro, Ho}, where it was also shown that the
fundamental group must be a free group. In contrast, in dimension $3$ any
curvature tensor can be realized locally, see \cite{BKV, Bry1}. But in
higher dimensions, the condition that $R$ needs to satisfy seems to be
quite strong. Surprisingly, the only known compact non-homogeneous
examples are the Ferus-Karcher-M\"unzner isoparametric hypersurfaces
(i.e., the eigenvalues of the shape operator are constant), see
\cite{FKM}.

It is natural to attack this kind of problem in the context of
submanifolds in space forms since, by the Gauss equation, more control is
gained on the structure of the curvature tensor. This is particularly
true for hypersurfaces where the curvature tensor, which in general is a
very complicated algebraic object, is determined just by an endomorphism
of the tangent bundle, namely, the shape operator.

In \cite{Ts} K. Tsukada studied the problem of classifying curvature
homogeneous hypersurfaces $M^n$ in the simply connected space form
$\Q^{n+1}_c$ of constant curvature $c$. He showed that any such
hypersurface is isoparametric, or has constant curvature~$c$, or has
rank two, that is, the rank of its shape operator is two everywhere; see
\pref{fbf}. Observe that constant scalar curvature and curvature
homogeneous are equivalent notions for rank two submanifolds in space
forms by the Gauss equation.

There are two natural examples of 3-dimensional complete rank two
hypersurfaces with constant scalar curvature,
$f_c:\Lambda\to\Q_c^4$, $c=\pm 1$. The first one $f_1$ is the unit normal
bundle $\Lambda$ of the Veronese surface $\RP_{1/3}\subset\Sp^4$,
which is one of Cartan's isoparametric hypersurfaces with three different
principal curvatures. This hypersurface is not only isoparametric, but
homogeneous. The second one $f_{-1}$ is the unit normal bundle $\Lambda$
of the flat torus in the De-Sitter space,
$g=(g_0,1)\colon T^2\to\DS^4\subset\R^{4,1}$, where\linebreak
$g_0\colon T^2\to\Sp^3(\sqrt{2})\subset\R^4\times\{0\}\subset\R^{4,1}$
is the minimal equivariant flat Clifford torus; see \sref{clif}.
This hypersurface has a two parameter family of symmetries induced
by the symmetries of the Clifford torus.

If $n\ge 4$ or $n=3$ and $c=0$ the condition on the hypersurface
is quite rigid. Tsukada showed that, in this case, besides the
obvious Euclidean cylinders over constant curvature surfaces, there is
only one rank two example, a complete hypersurface in the hyperbolic
5-space $\Hi^5$. We will give a simpler proof of this fact in \sref{h5}
together with a more geometric description of this example closely
related to both $f_c$'s.

On the other hand, the case $n=3$ and $c\neq 0$ remained an open problem,
see e.g. \cite{BKV} p.255 and \cite{CMP}.
Our purpose in this paper is to answer this question. Recall that a rank
two hypersurface in $\Q_c^4$ is foliated by special geodesics, the so
called relative nullity leaves, tangent to the kernel of the shape
operator.

\begin{theorem}\po\label{main}
Let ${\cal M}$ be the set of immersed rank two hypersurfaces in $\Q_c^4$,
$c=\pm 1$, whose induced metric has constant scalar curvature. Then
${\cal M}$ contains $f_c$ as the only complete example, an isolated
hypersurface $\hat f_c$ with a circle of symmetries, and a one parameter
family of hypersurfaces admitting no continuous symmetries.
Moreover, up to a covering, any connected hypersurface in ${\cal M}$
is an open subset of one of these, provided it has no leaf of relative
nullity of minimal points in the case $c=1$.
\end{theorem}

To prove \tref{main} we will make use of the Gauss Parametrization that
we recall in \sref{gps}, which
is a powerful tool to study hypersurfaces of constant rank in space
forms. Our hypersurfaces will then be the unit normal bundles of their
polar surfaces in $\DSc^4$, in fact globally since we will show that the
relative nullity geodesics are complete. These surfaces are characterized
by the property that all shape operators along unit normal directions
have the same non-zero determinant, and thus have constant Gaussian
curvature. Our work will then reduce to classifying such surfaces.
Topologically, the polar surfaces of the one parameter family of
hypersurfaces with no continuous symmetries in \tref{main} are
diffeomorphic to a pair of pants if $c=1$, or to either a cylinder or a
plane if $c=-1$; see \sref{es}. Moreover, the hypothesis on the minimal
points for $c=1$ is equivalent to asking for the polar surface to have no
minimal points.

This approach also allows us to get simple explicit
parametrizations for $\hat f_c$ in \tref{main} as follows. Set
$r_0:=\arccos(\sqrt{2/3})$, and $D^2(r_0)\subset \R^2$ the 2-disk of
radius $r_0$ with polar coordinates $(r,\theta)$ if $c=1$. Then, up to
congruences,
$\hat f_1 = \hat f_1(r,\theta,\alpha)\colon D^2(r_0)\times S^1\to\Sp^4$
and $\hat f_{-1} = \hat f_{-1}(r,\theta,\alpha)
\colon (r_0,\pi-r_0)\times S^1\times\R\to\Hi^4$
have the unified expression
\begin{equation}\label{fc1}
\hat f_c =
\begin{pmatrix}
\cos_c(\alpha)\sin(2\theta)\cos(2r)
- \sin_c(\alpha)\cos(2\theta)\cos(r)\sqrt{(3\cos(2r)-1)/2c}\\
\cos_c(\alpha)\cos(2\theta)\cos(2r)
+ \sin_c(\alpha)\sin(2\theta)\cos(r)\sqrt{(3\cos(2r)-1)/2c}\\
\cos_c(\alpha)\ \sin(\theta)\ \sin(2r)
- 2\sin_c(\alpha)\cos(\theta)\sin(r)\sqrt{(3\cos(2r)-1)/2c}\\
\cos_c(\alpha)\ \cos(\theta)\ \sin(2r)
+ 2\sin_c(\alpha)\sin(\theta)\sin(r)\sqrt{(3\cos(2r)-1)/2c}\\
(3/2)\sin_c(\alpha)(\cos(2r)-1)
\end{pmatrix},
\end{equation}
where $\sin_c$ and $\cos_c$ stand for $\sin$ and $\cos$ if $c=1$, $\sinh$
and $\cosh$ if $c=-1$. We will see that at the boundary $\hat f_1$
has a 2-torus as singular set (where its mean curvature is unbounded),
and at the boundary $\hat f_{-1}$ has two singular 2-cylinders.
In addition, from~\eqref{fc1} it is not hard to check that $\hat f_c$ is
algebraic; see \sref{rotsim}.

\smallskip

In a forthcoming paper \cite{Bry2}, a more complete description of the
corresponding polar surfaces will be provided. When $c=1$, it will be
shown that there exists a $1$-parameter family of real-analytic mappings
$\overline g_a:S^2\to S^4$ for $0\le a\le 1$ such that the polar surface
of a hypersurface as described in \tref{main} is congruent to an open
subset of $\overline g_a(S^2)$ for some $a$. The map $\overline g_a$ is a
topological embedding and is an immersion except along the equator
in~$S^2$, where its differential has rank~$1$. When $a=0$, the image
$\overline g_0(S^2)$ has a rotational symmetry and is congruent to the
algebraic surface described by \eqref{Y}. Its only minimal points are the
two `poles’ of the rotational symmetry. When $a>0$, the image
$\overline g_a(S^2)$ has an $8$-fold discrete symmetry group, and it
contains exactly four distinct minimal points (at which the surface
$\overline g_a(S^2)$ is smooth). At present, it is not known whether the
compact surface $\overline g_a(S^2)$ is algebraic when $a>0$. Meanwhile,
when $c=-1$, a correspondingly complete description will be given of the
polar surfaces of the hypersurfaces described by \tref{main}. Again, it
turns out that there is a $1$-parameter family of such polar surfaces up
to congruence, and, except for one particular value of the parameter, the
analytically-completed surfaces have similar singularity properties,
while, for the exceptional value, the singular structure is quite
different. Again, it is not known at present whether these surfaces are
algebraic.
\smallskip

The paper is organized as follows. In \sref{pre} we reduce the
classification to hypersurfaces of $\Q^{4}_c$ and explain the Gauss
parametrization. We then convert our problem to a classification of the
corresponding polar surface. In \sref{struc} we discuss the structure
equations of the polar surface and in \sref{compt} the compatibility
condition that needs to be satisfied in order to find local solutions. In
\sref{es} we determine the maximal domain of the polar surface, and its
topological type. Finally, in \sref{rotsim} we discuss the example with
rotational symmetry, and in \sref{h5} give a simple description of the
Tsukada example.

\section{Preliminaries}\label{pre} 
Let $M^n$ be a curvature homogeneous Riemannian manifold with curvature
tensor $R$. Let $f:M^n\to \Q^{n+p}_c$ be an isometric immersion
with second fundamental form $\alpha$ into the simply connected
space form $\Q^{n+p}_c$ of curvature $c$. Fix $x_0\in M^n$.
Then, for each $x\in M$, there is a
linear isometry $J_x:T_xM\to T_{x_0}M$ such that
\begin{equation}\label{ch}
R_x = J_x^* R_{x_0}.
\end{equation}
We say that such an $f$ is {\it weakly isoparametric} if for each
$x\in M$ there exists another linear isometry
$\hat J_x:T^\perp_xM\to T^\perp_{x_0}M$ such that
\begin{equation}\label{iso}
\hat J_x\circ\alpha_x = J_x^* \alpha_{x_0}.
\end{equation}
Notice that, by the Gauss equation, \eqref{iso} implies \eqref{ch}.

For each $x\in M$, define the bilinear map
$$
\beta_x:T_xM\times T_xM \to W_x:=T^\perp_xM\times T^\perp_{x_0}M
$$
as $\beta_x = (\alpha_x,J_x^* \alpha_{x_0})$.
Again by Gauss equation, $M^n$ is curvature homogeneous (with
respect to $J$) if and only if $\beta_x$ is {\it flat}, that is,
$$
\la \beta_x(X,Y),\beta_x(U,V)\ra = \la \beta_x(X,V),\beta_x(U,Y)\ra,
\ \ \ \forall\ X,Y,U,V\in T_xM,
$$
where the inner product on $W_x$ is the natural indefinite one of type
$(p,p)$, namely,
$\la\,,\,\ra=\la\,,\,\ra_{T^\perp_xM} - \la\,,\,\ra_{T^\perp_{x_0}M}$.
It turns out that $f$ is weakly isoparametric at $x$ if and only if
$\beta_x$ is {\it null}, i.e.,
$$
\la \beta_x(X,Y),\beta_x(U,V)\ra = 0,\ \ \ \forall\ X,Y,U,V\in T_xM.
$$
Indeed, for the converse just observe that the expression in \eqref{iso}
serves as a good definition of $\hat J_x$ between the images of
$\alpha_x$ and $J^*_x\alpha_{x_0}$, which can afterwards be extended
by linearity as a linear isometry.

Deciding when a flat bilinear map is null is a key point in isometric
rigidity problems of submanifolds. 
Theorem 3 in \cite{df} ensures that, if not null, a symmetric flat
bilinear form must have a highly degenerate component, at least if
$p\leq 5$. More precisely, $W_x$ decomposes orthogonally as
$$
W_x=W_0\oplus^\perp W_1
$$
and $\beta$ decomposes accordingly as $\beta = \beta_0 + \beta_1$,
where $\beta_0$ is null and $\beta_1$ has nullity of dimension
$\nu_x\geq n-\dim W_1$. In particular, if the codimension $p$ is equal
to 1 we have $\nu_x\geq n-2$.
We conclude the following (see Theorem 2.3 in \cite{Ts}):
\begin{proposition}\po\label{fbf}
A hypersurface in $\Q^{n+1}_c$
is curvature homogeneous if and only if it is isoparametric, or
has constant sectional curvature $c$, or has rank two with constant
scalar curvature.
\end{proposition}
\begin{remark}\po
{\rm The case of constant curvature $c$ is well understood, since the set
of such nowhere totally geodesic hypersurfaces can be naturally
parametrized by the set of regular smooth curves in $\Q^{n+1}_c$ using
the Gauss Parametrization; see \sref{gps}. The isoparametric case was
completely classified by E. Cartan for $c\leq 0$, while the $c>0$ case in
full generality still remains a well-known open problem.}
\end{remark}

In view of this, we concentrate from now on to the general task of
describing constant scalar curvature rank two hypersurfaces in space
forms $\Q_c^{n+1}$, of any dimension. Since the Euclidean case was solved
in \cite{Ts}, we restrict ourselves to the cases $c=\pm1$.

\subsection{The Gauss parametrization for rank two hypersurfaces}\label{gps} 
The Gauss parametrization is a powerful tool to work with hypersurfaces
with constant rank in space forms, as in our situation.
It was created by Sbrana in \cite{sb} with the purpose of classifying
nonflat locally isometrically deformable Euclidean hypersurfaces, which
also have constant rank two. The tool was studied in further detail in
\cite{dg}, and we briefly describe it next.

\medskip

Fix $c=\pm1$ and let $\E^n$ denote the corresponding Euclidean space
$\R^n$ or the Lorentzian space $\R^{n-1,1}$, i.e., $\R^n$ with the metric
$\d x_1^2+\cdots+\d x_{n-1}^2+c\, \d x_n^2$.
Let $f:M^n\to\Q^{n+1}_c\subset \E^{n+2}$, $n\geq 3$, be a rank two
connected orientable hypersurface and $\Delta^{n-2}$ its totally geodesic
relative nullity foliation, namely, the integral leaves of the kernel of
its second fundamental form. Consider the map
$\hat g:M^n\to\DSc^{n+1}:=\{x\in\E^{n+2}:\la x,x \ra=1\}$
such that $\{f,\hat g\}$ is an oriented pseudo-orthonormal
normal frame of $f$ seen in $\E^{n+2}$, namely, $\la\hat g, f\ra=0$ and
$\la\hat g(x), f_{*x}v\ra=0$ for all $x\in M^n, v\in T_xM$.
If we take the (local) leaf space
$$
\pi:M^n\to V^2:=M^n/\Delta,
$$
the map $\hat g$ descends to the quotient.
That is, there is an immersion called the {\it polar map} of $f$ given by
$$
g:V^2\to\DSc^{n+1}\ \ \ {\rm with}\ \ \ g\circ\pi=\hat g.
$$
We fix on $V^2$ the metric induced by $g$, which is
Riemannian since
\begin{equation}\label{par}
    \Delta^\perp(w)= g_{*p}(T_pV),\ \ \ p=\pi(w).
\end{equation}
It turns out that, locally, $f(M^n)$ can be seen as the unit normal
bundle $\Lambda$ of $g$,
$$
\Lambda:=\{w\in T^\perp_ g V\subset
T\DSc^{n+1}:p\in V^2,\ \la w,w \ra=c\,\},
$$
that is, as the image of the map $\hat f: \Lambda \to\Q^{n+1}_c$ that
sees each $w\in\Lambda$ as an element in $\Q^{n+1}_c\subset \E^{n+2}$
under parallel translation,
\begin{equation}\label{gp2}
\hat f(w) = w.
\end{equation}
The leaves of relative nullity of $f$ are then identified to (open
subsets of) the fibers of $\Lambda$ as a bundle. Denote by $A_w$ the
shape operator of $g$ in the direction
$w\in \Lambda_p\subset T^\perp_{g(p)}V$, $p=\pi(w)$.
It is easy to check that the
regular points of the {\it Gauss parametrization} $\hat f$ in \eqref{gp2}
are the vectors $w\in\Lambda$ such that $A_w$ is invertible.
Moreover, using the identification in \eqref{par}, we have that the shape
operator of $f$ restricted to $\Delta^\perp(w)$ in the direction
$\hat g$ is just $A_w^{-1}$. Conversely, given any surface
$g:V^2\to\DSc^{n+1}$, the map \eqref{gp2} gives a rank two hypersurface
in $\Q^{n+1}_c$, when restricted to the open subset of its regular points
as described above; see \cite{dg} for details.

\medskip

In view of this construction and the Gauss equation we are able
to transfer our problem to the polar map $g$:

\begin{proposition}\po\label{gp}
For $c=\pm1$, consider a connected orientable rank two hypersurface
$f:M^n\to\Q^{n+1}_c$ with polar map $g:V^2\to\DSc^{n+1}$. Then $f$ is
curvature homogeneous if and only if the map $w\in\Lambda\mapsto\det A_w$
is a non-zero constant.
\end{proposition}

Clearly, if this last map is constant along a small segment of a relative
nullity geodesic $\gamma$, then it must be constant along the whole
$\gamma$, and in particular $\hat f$ in \eqref{gp2} must be regular along
all of $\gamma$. We conclude that we can assume from now on that all
these geodesics are complete, even though we do not ask for the
hypersurface $M^n$ itself to be complete. In particular, $M^n$ becomes
the total space of the bundle $\Q^{n-1}_c\to M^n\xrightarrow{\pi} V^2$,
and we conclude the following.

\begin{corollary}\po\label{max}
For $f$ as above we have that $f(M^n)=\hat f(\Lambda)$. Conversely,
if a surface $g$ satisfies the property of the last proposition, then
$\hat f$ in \eqref{gp2} gives globally a curvature homogeneous
rank two immersion defined on the whole unit normal bundle $\Lambda$ of
$g$.
\end{corollary}
\proof
For each $p\in V^2$, since along an open subset of $\Lambda_p$ the map
$w\in\Lambda_p\mapsto\det A_w$ in \pref{gp} is a non-zero constant,
then this property holds over the whole leaf $\Lambda_p$, and therefore
the map $\hat f$ is an immersion over all of $\Lambda$, still
preserving the property that $\det A_w$ is a non-zero constant.
Therefore $\hat f$ extends the original immersion $f$ and is a curvature
homogeneous regular hypersurface defined over the whole bundle $\Lambda$.
\qed
\vspace{1.5ex}

As a consequence we conclude that our problem is low dimensional
(the case $n=4$ and $c=-1$ will be completely classified in \sref{h5}):

\begin{corollary}\po
Either $n=3$ for $c=1$, or $3\leq n \leq 4$ for $c=-1$.
\end{corollary}
\proof
Fix $p\in V^2$ and let $E^3$ be the vector space of self-adjoint
endomorphisms of $T_pV\cong\R^2$. A maximal subspace of
non-singular elements in $E^3$, excluding 0, has dimension~2.
Consider the linear map $\varphi: T^\perp_pV\to E^3$, $\varphi(w)=A_w$.
For $c=1$, clearly $\varphi$ has trivial kernel
since $\varphi(\Lambda_p)\subset {\rm Iso}(T_pV)$.
For $c=-1$, $\varphi$ has kernel at most one dimensional since
${\rm det}\circ\varphi$ is a non-zero constant (see e.g. \eqref{cor6}
in \sref{h5}). The proof follows from
$n=\dim\ker\varphi+\dim\im\varphi+1\leq \dim\ker\varphi+3$.
\qed

\begin{remark}\po
{\rm As already pointed out, the problem for hypersurfaces in Euclidean
space is simpler and completely understood. It turns out that the only
examples are $(n-2)$-cylinders over surfaces with constant Gaussian
curvature in $\R^3$, which are themselves classified. This can also be
easily obtained using the Gauss parametrization; see Theorem 3.4 in
\cite{dg}.}
\end{remark}

\begin{example}\po
{\rm
A well-known example of the situation in \pref{gp} is the minimal
Veronese surface $g_1:\RP^2_{1/3}\to\Sp^4\subset\R^5$.
It has the property that, given any orthonormal local tangent frame
$\{e_1,e_2\}$ of $T\RP^2_{1/3}$, there exists a unique orthonormal normal
frame $\{\xi_1,\xi_2\}$ such that
\begin{equation}\label{vl}
A_{\xi_1}=a\begin{pmatrix}0&1\\1&0\end{pmatrix},\ \
A_{\xi_2}=a\begin{pmatrix}1&0\\0&\!\!-1\end{pmatrix},\ \ \ a>0.
\end{equation}
A point on a surface in $\Sp^4$ whose second fundamental form satisfies
\eqref{vl} will be called {\it Veronese-like}. For $g_1$ we have that
$a=1/\sqrt{3}$, and $g_1$ is the only surface in $\Sp^4$ which is
Veronese-like everywhere. In fact, E. Cartan in \cite{ca} classified all
isoparametric hypersurfaces in space forms with 3 different principal
curvatures, the unit normal bundle of $g_1$ being the only one with
rank two.
}
\end{example}

The next lemma will be needed in the following sections. It is convenient
to call a basis $\{\xi_1,\xi_2\}$ of $\E^2$ orthonormal if
$\la\xi_i,\xi_j\ra=\delta_{ij}$ for $c=1$, as usual, while
$\epsilon:=\la\xi_1,\xi_1\ra=-\la\xi_2,\xi_2\ra=\pm 1$,
$\la\xi_1,\xi_2\ra=0$ if $c=-1$.

\begin{lemma}\label{al}\po
Let $g:V^2\to\DSc^4\subset\E^5$ be an isometric immersion such that
$\det A_w\neq 0$ is constant for all $w\in\Lambda$. Then, locally
around each non-minimal point of $g$, there exists an orthonormal
tangent frame $\{e_1,e_2\}$, an orthonormal normal frame
$\{\xi_1,\xi_2\}$, a constant $a>0$, and a smooth function $h>0$
on $V^2$, with $h>1$ if $c=1$, such that, in those frames,
\begin{equation}\label{al1}
A_{\xi_1}=a\begin{pmatrix} 0&1\\1&0 \end{pmatrix},\ \ 
A_{\xi_2}=a\begin{pmatrix}h&0\\0&\!\!-c/h\end{pmatrix}.
\end{equation}
Moreover, the Gaussian curvature of $V^2$ is constant $1-2\epsilon a^2$,
and outside the minimal points all this data is unique up to signs and
permutations of $e_1$ and $e_2$.
\end{lemma}
\proof
Let $L$ be the line bundle $L=\{\xi\in T^\perp_gV: \tr A_\xi=0\}$.
First, we claim that if $c=-1$ then $L$ is not light-like.
To see this,
assume otherwise, take a generator $\eta_1$ of $L$ and complete it to a
basis $\{\eta_1,\eta_2\}$ such that
$\la\eta_1,\eta_1\ra=\la\eta_2,\eta_2\ra=0$ and $\la\eta_1,\eta_2\ra=1$.
Then $\Lambda$ can be written as
$\Lambda=\{\eta_t=(t^{-1}\eta_1-t\eta_2)/\sqrt{2}: 0\neq t\in\R\}$. Write
$A_{\eta_1}=\tiny{a\!\begin{pmatrix}0&\!\!1\\1&\!\!0\end{pmatrix}}$
in some orthonormal basis, and
$A_{\eta_2}=\tiny{a\!\begin{pmatrix}x&\!\!y\\y&\!\!z\end{pmatrix}}$.
Thus, $2a^{-2}\det A_{\eta_t}=t^2(xz-y^2)+2y-1/t^2$,
which is not independent of $t$.

Now, choose $\xi_1\in L$ with $\la\xi_1,\xi_1\ra=\epsilon=\pm 1$
and complete it to an orthonormal normal frame $\{\xi_1,\xi_2\}$,
that is, $\la\xi_1,\xi_2\ra =0$ and $\la\xi_2,\xi_2\ra=\epsilon c$,
where of course $\epsilon=1$ if $c=1$. We can then write
$\Lambda=\{\xi_t=C_t\xi_1+S_t\xi_2: t\in I\subset \R\}$,
where $C_t$ and $S_t$ are smooth functions of $t$ satisfying
$c\,C_t^2+S_t^2=\epsilon$. In an orthonormal tangent frame of isotropic
vectors for $A_{\xi_1}$ we have that
$$
A_{\xi_1}=\begin{pmatrix} 0& a\\ a&0 \end{pmatrix},\ \ 
A_{\xi_2}=\begin{pmatrix}\alpha&\beta\\\beta&\gamma\end{pmatrix}.
$$
Hence, $\det A_{\xi_t}=(\gamma\alpha-\beta^2+c a^2)S_t^2
-2a\beta S_tC_t -c\epsilon a^2$, which must be a non--zero constant.
Therefore, $a\neq0$ is constant, $\beta=0$, and
$\gamma\alpha=-ca^2$. The lemma now follows easily.
\qed
\begin{remark}\po
{\rm
Notice that the set $\Sigma$ of minimal points of $g$ is nonempty only
if $c=1$ and it represents those points for which $h\to1$ in~\eqref{al1}.
Therefore all minimal points are Veronese-like. Hence, the problem
with the minimal points is that {\it any} pair of orthogonal tangent
directions provides the same normal form~\eqref{al1}, thus they are not
unique and the special frames in \lref{al} may not extend smoothly or
continuously to the minimal points, even if isolated.}
\end{remark}

\begin{remark}\po
{\rm
Since the shape operator of $f$ restricted to $\Delta^\perp$ at
the point $\xi_t\in\Lambda$ is $A_{\xi_t}^{-1}$, its mean curvature
is $c\epsilon S_t(h^2-c)/ha$. Thus, in terms of the Gauss parametrization,
the set of minimal points of (a maximal) hypersurface $f$, for $c=1$, is
$\Lambda|_\Sigma$ together with the two surfaces $\{\pm\xi_1(p): p\in
V^2\}\subset\Lambda$. In particular, $\Lambda|_\Sigma$ corresponds to the
set of leaves of relative nullity of $f$ contained in its set of minimal
points. Therefore, the exclusion of this set in \tref{main} is equivalent
to the exclusion of the minimal points of $g$.
}
\end{remark}

\begin{remark}\po\label{syms}
{\rm Since $V^2$ has constant Gaussian curvature it has many local
isometries. Yet, since $h$ and the frames in \lref{al} are unique, any
continuous family of (extrinsic) symmetries preserving $V^2$ cannot fix
points in $V^2$.}
\end{remark}

\section{Reduction of the structure equations}\label{struc} 
In this section we compute the structure equations of the polar map
$g:V^2\to\DSc^4\subset\E^5$ of our hypersurface in $\Q^4_c$, namely,
a nowhere minimal Riemannian surface as in \lref{al}.
For this section, and the remainder of the paper, the reader may find
it useful to verify our long but straightforward computations using
the Maple file \cite{BFZ}.

\vs

Following the notations in \lref{al}, extend the tangent frame
$\{e_1,e_2\}$ with $e_0:=g$, $e_3=\xi_1$, $e_4=\xi_2$. This is an
orthonormal frame of $\E^5$, since $\la e_i,e_j\ra=0$ if $i\neq j$ and
$$
\la e_0,e_0\ra=\la e_1,e_1\ra=\la e_2,e_2\ra=1,\ \ 
\la e_3,e_3\ra= c \la e_4,e_4\ra=\epsilon=\pm 1,
$$
with $c=\pm 1$, and $\epsilon=1$ if $c=1$.
Set
$$
\d e_i = \sum_j e_j \eta_{ji},\ \ \ \ {\rm with}\ \ \ 
\eta_{ji} \la e_i,e_i\ra = -\eta_{ij} \la e_j,e_j\ra,\ \ 0\leq i,j\leq 4.
$$
Hence $\eta_{30}=\eta_{40}=0$, and $\omega_1:=\eta_{10}$,
$\omega_2:=\eta_{20}$
must be linearly independent. The associated tangent and normal
connection 1-forms are $\omega:=\eta_{21}$ and $\mu:=\eta_{43}$, respectively.
\lref{al} is then equivalent to
$$
\eta_{13}=\epsilon a \,\omega_2,\ \eta_{14}=\epsilon c a h \,\omega_1,\ 
\eta_{23}=\epsilon a \,\omega_1,\ \eta_{24}=-\epsilon a h^{-1} \omega_2,
$$
with $a>0$ constant and $h>0$ smooth on $V$, with
$h>1$ if $c=1$ since we exclude minimal points.
Putting the above together gives
\begin{equation}\label{ee}
\eta=
\begin{pmatrix}
0&-\omega_1&-\omega_2&0&0\\
\,\omega_1&0&-\omega&-a\,\omega_2&-ah\,\omega_1\\
\,\omega_2&\omega&0&-a\,\omega_1&ach^{-1}\omega_2\\
0&\epsilon a \,\omega_2&\epsilon a \,\omega_1&0&-c\mu\\
0&\epsilon c a h\,\omega_1&-\epsilon a h^{-1}\omega_2&\mu&0
\end{pmatrix}.
\end{equation}

For convenience call $t_0=h$ and write
$$
\d t_0 = t_0(t_1\,\omega_1 + t_2\,\omega_2)
$$
for certain smooth functions $t_1,t_2$.
Recall that the structure equations are 
\begin{equation}\label{eseq}
\d \eta_{ji}=-\sum_k\eta_{jk}\wedge\eta_{ki}.
\end{equation}
These for $j=3,4$ and $i=1,2$ are the Codazzi equations. It is easy to
check using \eqref{ee} that they are equivalent to the determination
of the tangent and normal connections with the above data:
\begin{equation}\label{tnc}
\omega = -\frac{t_0^2t_2}{t_0^2-c}\,\omega_1+\frac{ct_1}{t_0^2-c}\,\omega_2,\ \ \ \ \
\mu = \frac{2ct_0^3t_2}{t_0^2-c}\,\omega_1-\frac{2ct_1}{t_0(t_0^2-c)}\,\omega_2. 
\end{equation}
For example, since $\eta_{31}=a\eta_{20}$, we get
$a \omega_1\wedge \omega-at_0 \mu\wedge \omega_1 = a \omega\wedge
\omega_1$, i.e., $\omega_1\wedge (2\omega+t_0\mu)=0$.

Now define the functions $t_{rs}$, $1\leq r,s\leq 2$ by
$$
\d t_i = t_{i1} \,\omega_1 + t_{i2} \,\omega_2.
$$
The structure equations \eqref{eseq} for $(j,i) = (1,2)$ and $(3,4)$
can be solved for $t_{11}$ and $t_{22}$ in terms of the others as
\begin{equation}\label{t11}
t_{11} = \frac{
c\epsilon a^2(5t_0^4-4ct_0^2-1)+2c(t_0^4(t_2^2-1)-2t_1^2)+2t_0^2(2t_1^2+1)
}{2(t_0^2-c)},
\end{equation}
\begin{equation}\label{t22}
t_{22} = \frac{
\epsilon a^2(5-4ct_0^2-t_0^4)
+2ct_0^2(2t_2^2+1)
-2(2t_0^4t_2^2-t_1^2+1)
}{2t_0^2(t_0^2-c)}.
\end{equation}
Since $0 = \d(\d(\log t_0)) = (t_1t_2 + t_{12}-t_{21})\,\omega_1\wedge \omega_2$,
we express $t_{12}$ and $t_{21}$ in terms of a new function $t_3$ as
$$
t_{12} = \frac{t_3}{t_0} - t_1t_2\frac{t_0^2}{t_0^2-c},\ \ \ \ \ \ 
t_{21} = \frac{t_3}{t_0} - t_1t_2\frac{c}{t_0^2-c}.
$$
Moreover, the identities $\d(\d t_1) = \d(\d t_2) = 0$ are equivalent to
$$
\d t_3 = \left(ct_0^3t_2(9\epsilon a^2-4)+6t_1t_3\right) \omega_1 +
\left(ct_1(9\epsilon a^2-4)-4t_0t_2t_3\right) \omega_2.
$$

\section{Compatibility analysis}\label{compt} 
At this point, we have, on $V^2$, two 1-forms $\omega_1$ and $\omega_2$
which satisfy
\begin{equation}\label{dw}
\d \omega_1 = -\frac{t_0^2t_2}{t_0^2-c} \,\omega_1 \wedge \,\omega_2\ \ \
{\rm and} \ \ \
\d \omega_2 = \frac{ct_1}{t_0^2-c} \,\omega_1 \wedge \,\omega_2,
\end{equation}
and four functions $t_0, t_1, t_2$, and $t_3$, whose exterior derivatives
are expressed explicitly in terms of $\omega_1, \omega_2$ and
$t_0, t_1, t_2$, and $t_3$. In addition, it is easy to check that the
structure equations are satisfied by our choices.

By a theorem of \'Elie Cartan \cite{ca0}, if these explicit formulae
for the exterior derivatives imply that $\d(\d t_k) = 0$ for $k=0,1,2,3$,
then local solutions exist in the following sense:
for every set of constants $r = (r_0,r_1,r_2,r_3)$, with $r_0\geq1$
and $r_0>1$ if $c=1$, there exist a surface $V_r$ and a point
$p_r \in V_r$ such that, on $V_r$, there exist a
coframing $\omega_1, \omega_2$ and smooth functions $t_0, t_1, t_2$, and
$t_3$ such that $t_k(p_r) = r_k$. Moreover,
such a surface $V_r$ is unique up to local diffeomorphism fixing $p_r$.
Thus, if the $\d^2 = 0$ identity were to hold formally for this system,
there would be a 4-parameter family of germs of `solution manifolds'
to these differential equations. However, it turns out that $\d^2 = 0$
is not an identity for this system.

Of course, we know that we must have $\d(\d \omega_1)=\d(\d\omega_2) = 0$
and $\d(\d t_0) = \d(\d t_1) = \d(\d t_2) = 0$,
because we used those equations to find the formula for $t_{ij}$,
but we have not checked whether $\d(\d t_3)$ vanishes. In fact, it
turns out that the above formulae imply
$$
\d(\d t_3)=-\frac{R_0[a,t_0,t_1,t_2,t_3]}{2ct_0}\ \omega_1\wedge\omega_2,
$$
where
$$
R_0[a,t_0,t_1,t_2,t_3]=20ct_3^2-
(9\epsilon a^2-4)(\epsilon a^2(t_0^4+10ct_0^2+1)
-12(t_0^4t_2^2 + t_1^2)-4ct_0^2).
$$
Notice that, if $\epsilon=1$ and $a=2/3$, then the vanishing of $R_0$
is equivalent to the vanishing of $t_3$.
Consequently, we will obtain a system satisfying Cartan's
Conditions when $\epsilon=1$ and $a = 2/3$ by setting $t_3 = 0$.
We have shown:
\begin{proposition}\label{hay}\po
If $\epsilon=1$ and $a=2/3$ there exists precisely a 3-parameter family
of germs of non-minimal surfaces $g$ as in \lref{al} for both $c=1$ and
$c=-1$.
\end{proposition}

We now rule out the remaining cases.
\begin{proposition}\label{nohay}\po
Let $g$ be a non-minimal surface as in \lref{al}.
If either $\epsilon=-1$ or $a\neq 2/3$, then $c=-1$ and $h\equiv1$
is constant.
\end{proposition}
\proof
Let $R_0$ be the polynomial in $a, t_0, \dots , t_3$ defined above. This
polynomial vanishes on every solution to the structure equations, and
hence its exterior derivative does as well. Compute $\d(R_0)$ using the
formulae for the derivatives of the $t_k$. This will be a 1-form that is
a linear combination of $\omega_1$ and $\omega_2$ with coefficients that
are rational functions of $a, t_0, \dots , t_3$ with denominators that
are products of powers of $t_0$ and $t_0^2-c$. Let $R_1$ be the numerator
of the coefficient of $\omega_1$ in $\d(R_0)$ and let $R_2$ be the
numerator of the coefficient of $\omega_2$ in $\d(R_0)$. Then $R_1$ and
$R_2$ are polynomials in $a, t_0, \dots , t_3$ that vanish on all
solutions of the structure equations.

Continuing, let $R_{11}$ be the numerator of the coefficient of
$\omega_1$ in $\d(R_1)$ and let $R_{12}$ be the numerator of the
coefficient of $\omega_2$ in $\d(R_1)$, when these coefficients are
expressed as rational functions of $a, t_0, \dots , t_3$ with
denominators that are products of powers of $t_0$ and $t_0^2-c$.

In this way, we generate a sequence of polynomials
$R_0, R_1, R_2, R_{11},\dots$. Consider the ideal $F$ in the polynomial
ring $\R[a, t_0, \dots , t_3]$ generated by the 15 polynomials
$$
R_0, R_1, R_2, R_{11}, R_{12}, R_{21}, R_{22}, R_{111}, \dots , R_{222}.
$$
Let $B$ be the Groebner basis of this ideal computed using the pure
lexicographical order $t_3 > t_2 > t_1 > t_0 > a$. Then $B$ is
an ordered list with 39 elements. The fourth element of $B$ factors as
$$
B_4 = (t_0^2-c)(9\epsilon a^2 - 4)^2 P (a, t_0),
$$
where $P(a,t_0)$ is an irreducible polynomial of degree 16 in $a$ and
$t_0$ (the reader can verify this claim using \cite{BFZ}).
Now, $B_4$ being in the ideal $F$ must vanish on any solution of
the structure equations. Since $t_0^2 \neq c$, it follows that either
$a=2/3$ and $\epsilon=1$, or else $P(a,t_0)=0$.

However, if $P(a,t_0)$ vanishes identically on the solution, then $t_0$
must be a root of a nontrivial polynomial with constant coefficients and
hence $t_0$ must be constant. Since $\d t_0$ would then vanish identically,
it would then follow that $t_1$ and $t_2$, and hence $t_{11}, t_{12},
t_{21}$ and $t_{22}$ would vanish identically, but this is clearly
impossible unless $t_0\equiv 1$ by \eqref{t11} and \eqref{t22}.
\qed

\begin{corollary}\po
The Ricci eigenvalues of the hypersurfaces in \tref{main} are
constant. For $f_1$ they are $\{2,-1,-1\}$, for $f_{-1}$ they are
$\{-2,-4,-4\}$, while for all the others they are
$\{2c, 2c-9/4,2c-9/4\}$.
\end{corollary}
\proof
We only check for the last case, since $f_1$ is well-known,
and the principal curvatures of $f_{-1}$ will be computed in the next
section.

Since our 3 manifold has relative nullity $\Delta$ of dimension 1,
its corresponding Ricci eigenvalue is $2c$. Now, any vector in
$\Delta^\perp$ is then a Ricci eigenvector with eigenvalue equal to 
$c+(c-a^{-2})=2c-9/4$.
\qed

\subsection{The case \texorpdfstring{$c=-1$}{c=-1} and \texorpdfstring{$h\equiv1$}{h=1}}\label{clif} 

In this case, we have that $t_0=h\equiv1$, and then $t_i=t_{ij}=t_3=0$
and, by \eqref{dw}, $\d \omega_i=0$ for $i,j=1,2$. In addition,
$\omega=\mu=0$, so $V^2$ is a flat surface with flat normal bundle. In
particular, $\epsilon=1$ and $a=1/\sqrt{2}$ by the second part of
\lref{al}. Since
$\d\eta=-\eta\wedge\eta$, we conclude from Maurer-Cartan Fundamental
Lemma that there exists a unique (up to left translations) solution
$G:\tilde V^2\to \SO(4,1)$ of the system $\d G=G\eta$ defined on the
universal cover $\tilde V^2$ of $V^2$. In our situation, this is just
$G=e^\gamma$, where $\eta=\d\gamma$ and
$$
\gamma(x,y)=
\begin{pmatrix}
0 & -\sqrt{2} & 0 & 0 & 0\\
\sqrt{2} & 0 & 0 & 0 & -1 \\
0 & 0 & 0 & -1 & 0\\
0 & 0 & 1 & 0 & 0\\
0 &-1 & 0 & 0 & 0
\end{pmatrix}\frac{x}{\sqrt{2}}+ 
\begin{pmatrix}
0 & 0 & -\sqrt{2} & 0 & 0\\
0 & 0 & 0 & -1 & 0\\
\sqrt{2} & 0 & 0 & 0 & -1\\
0 & 1 & 0 & 0 & 0\\
0 & 0 & -1 & 0 & 0
\end{pmatrix}\frac{y}{\sqrt{2}}.
$$
Then $g=e_0(G)$ is the flat two torus
$g:T^2:=\R^2/\Z^2\to \DS^4\subset\R^{4,1}$ given by
$$
g(x,y):=
\begin{pmatrix}
2\cos(x)\cos(y)-1\\
-\sqrt{2}\sin(x)\cos(y)\\
-\sqrt{2}\cos(x)\sin(y)\\
\sqrt{2}\sin(x)\sin(y)\\
-\sqrt{2}\cos(x)\cos(y)+\sqrt{2}
\end{pmatrix},
$$
whose induced metric is twice the canonical one.
It is easy to check that $g$ satisfies \lref{al} with $h\equiv 1$
and $a=1/\sqrt{2}$.

Observe now that $g$ is also contained in the hyperplane
$x_1+\sqrt{2} x_5=1$. In fact, after a change of orthonormal basis $g$
can be written as
\begin{equation}\label{ct}
g=(g_0,1)\colon T^2\to\Sp^3(\sqrt{2})\times\R\subset\R^4\times\R=\R^{4,1},
\end{equation}
where 
$g_0\colon T^2\to\Sp^3(\sqrt{2})\subset\R^4$ is the standard minimal
equivariant flat Clifford torus,
$$g_0(x,y)=\sqrt{2}\,(\cos(x)\cos(y),\cos(x)\sin(y),\sin(x)\cos(y),\sin(x)\sin(y)).$$
Then, $g$ is an ${\rm Iso}(T^2)$-equivariant
isoparametric surface in codimension two in the De-Sitter space
$\DS^4$, i.e., it has parallel second fundamental form.
The corresponding non-isoparametric complete curvature
homogeneous hypersurface
$f_{-1}:\Lambda = T^2\times\R \to \Hi^4$ is thus given by
\begin{equation}\label{f1}
f_{-1}(x,y,t)= \sinh(t)\xi_1 + \cosh(t)\xi_2
=\frac{1}{\sqrt{2}}\left(\cosh(t)g_0+\sqrt{2}\sinh(t)\xi,2\cosh(t)\right),
\end{equation}
where $\xi=(g_0)_{xy}/\sqrt{2}$ is the Gauss map of $g_0$.
The equivariant isometries of $g_0$ induce a two-parameter family of
extrinsic symmetries of $f_{-1}$. The principal curvatures of $f_{-1}$
are $\{\cosh(t)+\sinh(t),\cosh(t)-\sinh(t),0\}$.

\section{Existence of solutions}\label{es} 
In this section we compute the maximal surfaces in \pref{hay}.

\vsss

As already seen, in this case we must have $\epsilon=1$, $a=2/3$ and
$t_3=0$, and therefore our system becomes
$$
\begin{pmatrix}
\d t_0\\ \d t_1\\ \d t_2
\end{pmatrix}=
\begin{pmatrix}
t_0t_1&t_0t_2\\
\frac{2(t_0^2-c)(9t_1^2+1)+ct_0^4(9t_2^2+1)-t_0^2}{9(t_0^2-c)}&
-\frac{t_0^2t_1t_2}{t_0^2-c}\\
-\frac{ct_1t_2}{t_0^2-c}&
\frac{-2t_0^2(t_0^2-c)(9t_2^2+1)+9t_1^2+1-ct_0^2}{9t_0^2(t_0^2-c)}\\
\end{pmatrix}
\begin{pmatrix}
\omega_1\\\omega_2
\end{pmatrix},
$$
together with the given formulae for $\d \omega_1$ and $\d \omega_2$ in
\eqref{dw}. Now, as one can verify, one has the identity $\d(\d t_k) = 0$
for $k = 0, 1, 2$, so Cartan's Theorem suffices to prove existence
of a one parameter family of surfaces since two degrees of freedom come
from moving the base point over the surface.

\vsss

Now, one can, in this case, prove existence without having to quote
Cartan's Theorem, at the price of doing some further computation. In fact,
there are other advantages to doing an explicit computation, as will be
seen.

\vsss

Let us write the above equation in the form
$$
t_0^2(t_0^2-c)(\d t_0,\d t_1,\d t_2)=
P[t_0,t_1,t_2]\ \omega_1+Q[t_0,t_1,t_2]\ \omega_2,
$$
where $P[t_0,t_1,t_2]$ and $Q[t_0,t_1,t_2]$ are $\R^3$-valued polynomials
in $t_0, t_1, t_2$. Define the \mbox{$\R^3$-valued} polynomial
$N[u_0, u_1, u_2]$ as
$$
t_0^2(t_0^2-c)N[t_0, t_1, t_2]=
P[t_0,t_1,t_2] \times Q[t_0,t_1,t_2].
$$
Notice that the entries of $N$ have no common factor by the definitions
of $P$ and $Q$.

Consider the 1-form $\theta$ on $\R^3$ defined by
$$
\theta = \la N [u], \d u\ra.
$$
where $[u]=[u_0, u_1, u_2]$ and $\d u=(\d u_0,\d u_1, \d u_2)$.
Calculation shows that $\theta$ vanishes only along the two curves
$$
C_1=\{u_2=0, 9u_1^2=(2u_0^{\ 2}+c)(u_0^{\ 2}-c)\},
$$
$$
C_2=\{u_1=0, 9u_2^2=(2u_0^{-2}+c)(u_0^{-2}-c)\},
$$
and these two curves only intersect when $c=1$ and do so at the points
$(u_0,u_1,u_2)=(\pm 1,0,0)$. Moreover, one computes that
$\theta\wedge \d \theta = 0$, i.e., the distribution
$$D=\ker \theta$$
on $\RR=\{u\in\R^3:u_0>0\}$ satisfies Frobenius integrability, so
that its leaves foliate $\RR \setminus (C_1 \cup C_2)$. In fact, a
calculation allows one to find a first integral. Indeed, setting
\begin{equation}\label{l}
L := \frac{u_0^4(u_0^2(9u_2^2+1)+c(9u_1^2+1))^2}
{(u_0^4(9u_2^2+1)+cu_0^2+(9u_1^2+1))^3}
\end{equation}
one gets that $\theta\wedge \d L = 0$. Note that $0\leq L \leq 4/27$, with
$L=4/27$ only on $C_1\cup C_2$. Moreover, $L=0$ only
when $c=-1$ and on the hypersurface
$\Omega=\{u_0^2=(9u_1^2+1)/(9u_2^2+1)\}\subset\R^3_+$, which is
homeomorphic to a plane. For any other value $0<R<4/27$, $L^{-1}(R)$ is a
smooth integral surface of $D$ which cannot intersect the plane $u_0=0$.

Notice also that $L$ is invariant under the transformation
$$
\varphi(u_0,u_1,u_2) =(1/u_0,u_2,u_1),
$$
and that $\varphi$ interchanges $C_1$ and $C_2$.
This corresponds to an arbitrary choice between $h\geq1$ and $h\leq 1$,
and the corresponding swap of the elements of the tangent frame in
\lref{al}.

For $c=1$, let $\Pi\subset\RR$ be the plane $u_0=1$ and $\Sigma\subset V$
the set of minimal points of~$g$. For $c=-1$, set both sets $\Pi$ and
$\Sigma$ as empty.

If $c=1$, all 2-dimensional leaves of $D$ intersect $\Pi$ transversally
since $\theta$ is nonvanishing when pulling back to $\Pi$. In fact,
given $r\geq0$, if $R:=(9r^2{+}2)^2/(9r^2{+}3)^3$ the intersection
$\Pi\cap L^{-1}(R)$ is the circle ${\cal C}_r$ of radius
$r$ centered at the origin, with $r\to0$ as $R\to4/27$ and
$r\to +\infty$ as $R\to 0$.
Each 2-dimensional leaf of $D$ is a union of 2 pair of pants
glued at their `waistline' ${\cal C}_r$ (rotated $90^\circ$
from being aligned with the 'legs' of the opposite pair), that are
interchanged by $\varphi$, and which then becomes a tube over the
connected curve $C_1\cup C_2$; see the picture on the left in Figure 1.

If $c=-1$, each 2-dimensional leaf $L^{-1}(R)$ for $0<R<4/27$ has two
connected components separated by $\Omega$, each of which is
diffeomorphic to a cylinder as a tube around one of the disjoint curves
$C_1$ or $C_2$; see the picture on the right in Figure 1.

\tagged{pics}{
\begin{figure}[ht!]
\centering{\includegraphics[width=75mm]{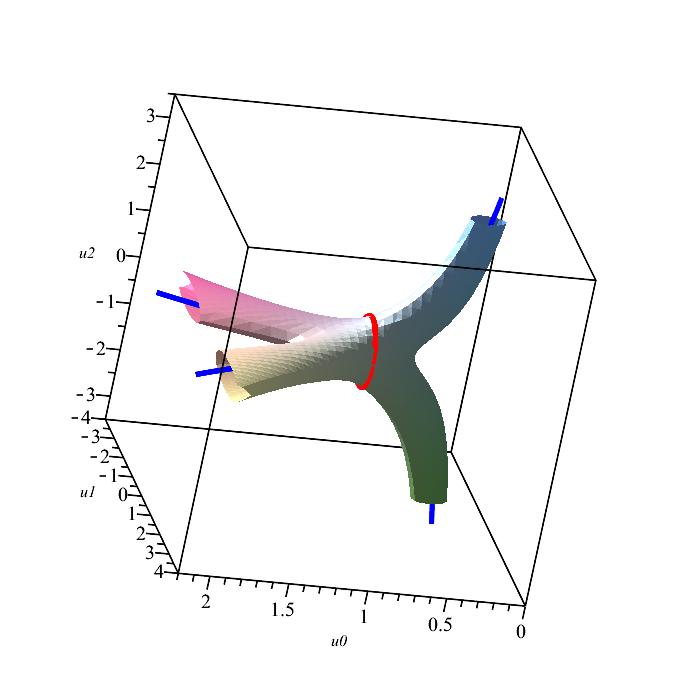}
\includegraphics[width=75mm]{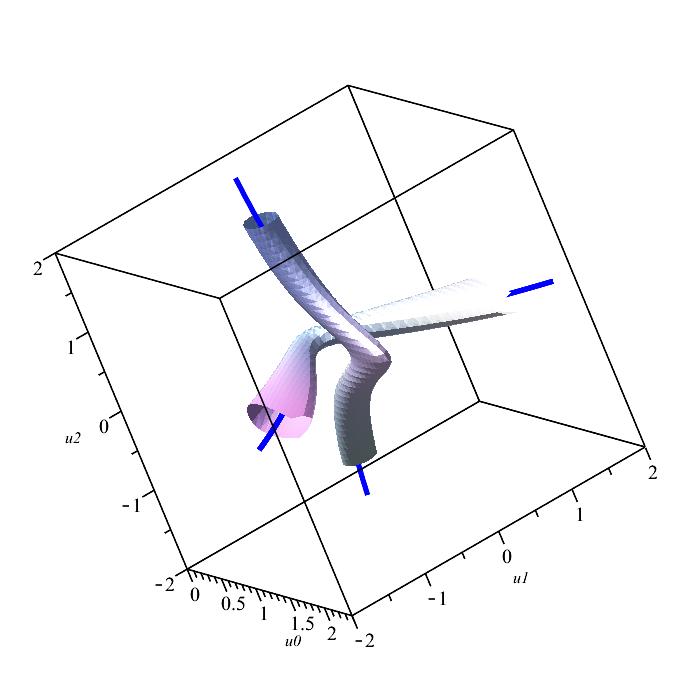}}
\caption{Leaves of the foliations $D$ for $c=1$ and $c=-1$}
\end{figure}
}

\vspace{.3cm}

Let $V^*$ be a connected component of $V\setminus\Sigma$.
By construction, since $N$ is perpendicular to both $P$ and $Q$,
the function $t=(t_0,t_1,t_2): V^*\to\RR\setminus\Pi$ pulls back $\theta$ to
zero, i.e., it maps $V^*$ onto a leaf
of $D$. Because $N[u_0, u_1, u_2]$ does not vanish outside $C_1 \cup C_2$,
it follows that the map $t:V^* \to\RR\setminus\Pi$ is an immersion unless
its image lies in either $C_1$ or $C_2$. Thus, unless
$t(V^*) \subset C_1 \cup C_2$, one can regard $V^*$, up to a covering,
as an open set in a leaf of $D$.

Conversely, if $V \subset \RR$ is a 2-dimensional leaf of $D$ then $u_0,
u_1$, and $u_2$ restricted to a connected component $V^*$ of
$V\setminus\Sigma$ define functions $0<t_0$, $t_1$ and $t_2$, with
$t_0\neq 1$ if $c=1$, such that the differential of $t=(t_0,t_1,t_2)$
satisfies $\la N[t_0,t_1,t_2], \d t\ra=0$. It follows that there will be
unique 1-forms $\psi_1$ and $\psi_2$ on $V^*$ satisfying $\d t=P[t]
\psi_1+Q[t] \psi_2$. Setting $\omega_i=t_0^2(t_0^2-c) \psi_i$, $i=1,2$,
then defines a coframe on $V^*$. One can verify that this coframe
satisfies \eqref{dw}.
In particular, now defining the various $\eta_{ab}$, $0 \leq a, b\leq 4$
using their formulae given above in terms of the $\omega_i$ and
$t_0,t_1,t_2$, the 1-form $\eta$ satisfies $\d\eta = -\eta\wedge\eta$. By
Maurer-Cartan Fundamental Lemma, there will be a mapping $G$ from the
simply connected cover $\tilde V^*$ of $V^*$ into $\SO_c(5)$, where
$\SO_c(5)=\SO(5)$ if $c=1$, or $\SO(4,1)$ if $c=-1$,
such that $G^{-1}\d G = \eta$. The resulting mapping
$g=e_0(G)\colon\tilde V^*\to\DSc^4$ will then give an immersion of
$\tilde V^*$ onto $\DSc^4$ as a surface satisfying \lref{al} with $a =
2/3$ and $\epsilon=1$.
Notice that, by the above discussion, $\tilde V^*$ is homeomorphic to the
universal cover of a pair of pants if $c=1$, and to a plane if $c=-1$.

Since there is a 1-parameter family of 2-dimensional leaves of $D$, these
give a \mbox{1-parameter} family of these surfaces in $\DSc^4$ that
have no continuous symmetries (since the map $t$ is an immersion and it
should be invariant by all symmetries, see \rref{syms}), and every such
connected surface in $\DSc^4$ without continuous symmetries is, locally,
an open set in one of these surfaces.

\medskip

It turns out that none of these surfaces is complete:

\begin{proposition}\label{noncomplete}\po
The only complete surfaces $g$ as in \lref{al} are the Veronese surface
and the torus in \eqref{ct}. In particular, there is no rank two
complete curvature homogeneous hypersurface in $\Q_c^4$ besides $f_c$.
\end{proposition}
\proof
Assume such a complete surface $g$ different from the Veronese and the
one in \eqref{ct} exists. Since $a=2/3$ and $\epsilon=1$, the Gaussian
curvature of $g$ is constant $1/9>0$. Hence the surface is diffeomorphic
to either $\Sp^2$ or $\RP^2$.

For $c=-1$, since there are no minimal points in $g$ we have a
global coframe $\omega_1,\omega_2$ on $\Sp^2$ which is obviously impossible.

For $c=1$, a computation shows that the square of the mean curvature 
vector of $g$, namely, $H=(h^2-1)^2/h^2$, is a superharmonic function,
since
$$
h^4\Delta H/2 = (4h^6+h^4+2h^2+1)t_1^2
+h^2(h^6+2h^4+h^2+4)t_2^2+(h^4-1)^2/9\geq 0.
$$
Thus $H$ and $h$ are constant. By the above $h=1$, $g$ is minimal and
therefore the Veronese surface.
\qed

\begin{remark}\po\label{for}
{\rm
In \cite{Bry2} the global topology of these surfaces will be addressed.
In particular, for $c=1$, it will be shown that $\Sigma$ is a smooth
isolated minimal point in $V$ and that the structure equations can be
extended smoothly to the circles ${\cal C}_r$.
}
\end{remark}

\section{The rotationally symmetric case}\label{rotsim} 
In this section we analyze the remaining case, namely, when $t(V)$
lies in one of the curves $C_1$ and $C_2$.

\vs

We first claim that we may assume that $t(V)\subset C_1$. Indeed, for
$c=-1$, the case $t(V)\subset C_2$ is completely analogous, since it
corresponds to reversing the roles between $e_1$ and $e_2$ (and thus
between $h$ and $-c/h$) in \lref{al}, namely, the $\varphi$-invariance
above. In particular, both cases give isometric surfaces. For $c=1$, the
curve $C_2$ is empty if $t_0>1$ by the sign of the right hand side
polynomial defining the curves.

Now, since $t_2=0$ and $t_1^2 = (2t_0^2+c)(t_0^2-c)/9$,
the structure equations are
$$
\d \omega_1=0,\ \ \d \omega_2 = \frac{ct_1}{t_0^2-c} \,\omega_1 \wedge \omega_2,\ \
\d t_0=t_0t_1\,\omega_1,\ \ \d t_1=\frac{1}{9}t_0^2(4t_0^2-c)\,\omega_1.
$$
Notice that the last one is a consequence of the third one and the
above formula for $t_1^2$. These can be easily solved for certain
coordinates $r$ and $\theta$ on $V$ as
$$
\omega_1=3\d r, \ \ \omega_2=3\sin(r) \d\theta,\ \
t_0=\sqrt{\frac{2c}{3\cos(2r)-1}},\ \ 
t_1=\frac{\sin(2r)}{3\cos(2r)-1}.
$$
Set $r_0=\arccos(\sqrt{2/3}\,)$ and $r_1=\pi-r_0$. A maximal domain of
the chart is $0<r<r_0$ if $c=1$ and $r_0<r<r_1$ if $c=-1$, namely,
$V=D^2(r_0)$ is a disk of radius $r_0$ if $c=1$ and the annulus
$V=(r_0,r_1)\times S^1$  if $c=-1$.
Moreover, the surface becomes singular as $r\to r_i$ where its mean
curvature vector field is unbounded. If $c=1$, then $r\to0$ if and only
if $t_0\to 1$, that is, the origin is the only minimal point of $V$, and
one can verify that it is a smooth point.

Using these formulae in \eqref{ee} and \eqref{tnc} we get
$\eta = \eta_1(r)\d r+\eta_2(r)\d\theta$, where
$$
\eta_1(r)=
\begin{pmatrix}
0 & -3 & 0 & 0 & 0\\
3 & 0 & 0 & 0 & -2 t_0 \\
0 & 0 & 0 & -2 & 0\\
0 & 0 & 2 & 0 & 0\\
0 & 2 ct_0  & 0 & 0 & 0
\end{pmatrix},\ \ 
\eta_2(r)=
\begin{pmatrix}
0 & 0 & -3\sinc & 0 & 0\\
0 & 0 & -\cosc & -2\sinc & 0\\
3\sinc & \cosc & 0 & 0 & 2c\sinc/t_0\\
0 & 2\sinc & 0 & 0 & 2c\cosc/t_0\\
0 & 0 & -2\sinc/t_0 & -2\cosc/t_0 & 0
\end{pmatrix},
$$
where $\sinc$ and $\cosc$ stand for $\sin(r)$ and $\cos(r)$ for
clarity.
It is easy to verify that the structure equation
$\d\eta=-\eta\wedge\eta$, or equivalently $[\eta_1,\eta_2]=-\eta_2'$, is
satisfied. Maurer-Cartan Fundamental Lemma thus implies that there is a map
$G_c:V\to \SO_c(5)$ such that $G_c^{-1}\d G_c = \eta$.
Then $\hat g_c=e_0(G_c): V \to \DSc^4$ gives
an immersion whose image is a surface in $\DSc^4$ as in \lref{al}.
Observe that $\hat g_c$ has a 1-parameter symmetry group induced by
translations in $\theta$, since $\eta$ is invariant under them. Notice
also that the system $G_c^{-1}\d G_c = \eta$ is equivalent to
$$
\frac{\partial G_c}{\partial r}=G_c\,\eta_1(r),\ \ \ 
\frac{\partial G_c}{\partial \theta}=G_c\,\eta_2(r). 
$$
The first equation (or equivariance) implies that
$G_c(r,\theta)=e^{\theta H}T(r)$ and $T'=T\eta_1$, with
$H\in \mathfrak{so}_c(5)$. By the second equation,
$H=T(r)\eta_2(r)T(r)^{-1}$ does not depend on $r$ and gives us $H$.
In addition, since $\eta_1=\eta_{11}\oplus\eta_{12}$ is
reducible in the $\{e_0,e_1,e_4\}$ and $\{e_2,e_3\}$ subspaces, the
problem becomes an ODE in $Gl(5,\R)$ of the form $T_1'=T_1\eta_{11}$
by taking an initial value in $\SO_c(5)$.
This is easily and explicitly integrable, giving $G$ whose first column
is (congruent to) the surface
\begin{equation}\label{gc}
\hat g_c = 
\begin{pmatrix}
3\sin(\theta)\sin(r)\cos(2r)\\
3\cos(\theta)\sin(r)\cos(2r)\\
(3/2)\sin(2\theta)\sin(r)\sin(2r)\\
(3/2)\cos(2\theta)\sin(r)\sin(2r)\\
\left((3\cos(2r)-1)/2c\right)^{3/2}
\end{pmatrix}.
\end{equation}

As a subset in $\R^5$, $\hat g_c(V)$ is cut out by three polynomial
equations, so it is contained in a singular algebraic surface
${\cal V}_c$, the intersection of three polynomials of degrees 2,~3,
and~6. The parametrization above only gives half of ${\cal V}_c$, the
half for which the fifth coordinate is greater than or equal to zero. The
other half is got by replacing the fifth coordinate with its negative.
The values for which $r\to r_i$ are closed torus knots in the Clifford
torus in~$\Sp^3$ where ${\cal V}_c$ has `creases', and it is smooth
everywhere else. It is here that the mean curvature of $\hat g_c$ goes to
infinity, similarly to the rim of the tractroid in $\R^3$ with
constant curvature~$-1$. The points $(0,0,0, 0,\pm1)$ for $c=1$ are the
minimal smooth points. Therefore, the maximal $V$ is a disk for $c=1$,
and an annulus for $c=-1$, with  torus knots as boundaries.
Moreover, making the substitution $u=\sqrt{3}\sin(\theta)\sin(r)$,
$v=\sqrt{3}\cos(\theta)\sin(r)$ and $w=\pm\sqrt{(3\cos(2r)-1)/2}$ in
\eqref{gc}, we see that for $c=1$ the full ${\cal V}_1$ is smoothly
parametrized by the unit sphere $u^2+v^2+w^2=1$ in the form
\begin{equation}\label{Y}
Y(u,v,w) = \frac{1}{\sqrt{3}}
\begin{pmatrix}
u\,(1+2w^2)\\
v\,(1+2w^2)\\
2uv\,\sqrt{2+w^2}\\
(u^2\!-\!v^2)\,\sqrt{2+w^2}\\
\sqrt{3}\,w^3
\end{pmatrix}.
\end{equation}
The embedding $Y$ is a smooth immersion away from the circle $w=0$,
which corresponds to the crease.

Notice also that, by \eqref{gc}, $\hat g_c$ can also be constructed as a
specific $\Sp^1$-orbit in $\R^4$ of a piece of the algebraic plane curve
$$
(x^2+4y^2)^3-9(x^2+4y^2)^2+81y^2=0,
$$ parametrized by
$\beta(r)=3\sin(r)(\cos(2r),\sin(2r)/2)$, and then simply adding as a
fifth coordinate $\sqrt{c(1-\|\beta\|^2)}$ to place it in $\Sp^4_c$. From
this we easily see that $\hat g_c$ is embedded; see Figure~2.

\tagged{pics}{
\begin{figure}[ht!]
\centering{
\includegraphics[width=60mm]{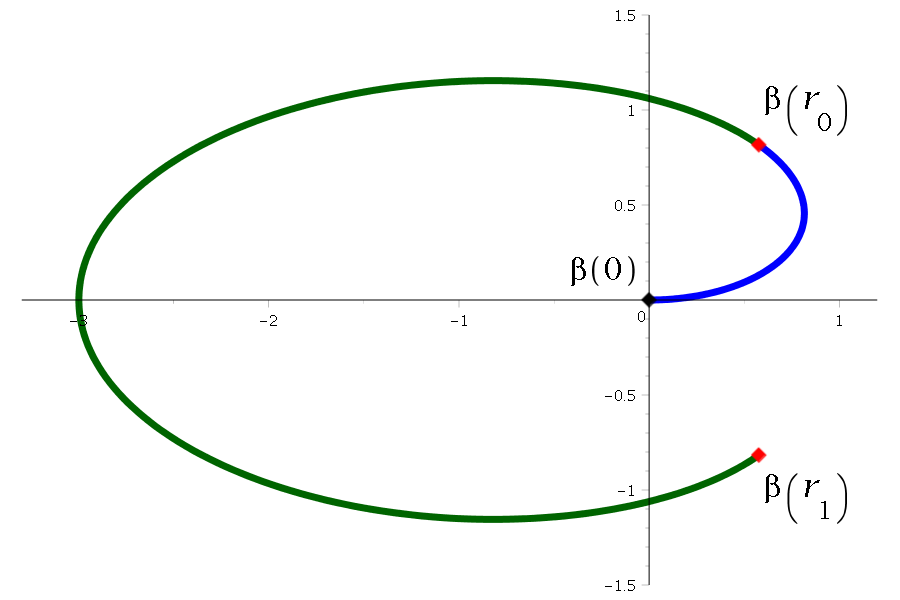}}
\caption{The curve $\beta$, in blue for $c=1$ and green for $c=-1$}
\end{figure}
}

\vspace{0.1in}

We can also use the map $G_c$ to give an explicit parametrization of
$\hat f_c$ in \tref{main} by taking
$\hat f_c=\hat f_c(r,\theta,\alpha)
=\cos_c(\alpha)e_3(G)+\sin_c(\alpha)e_4(G)$,
giving us the explicit expression \eqref{fc1} in the Introduction.
The image of $\hat f_c$ is also contained in an algebraic hypersurface
of $\Q^4_c$, namely, the intersection of $\Q^4_c\subset \E^5$ with the
0-level set of the polynomial
$$
64x_5^4(R+1)-\left(x_5^2(R^2-4R-8)-27c\,(x_1(x_3^2-x_4^2)+2x_2x_3x_4)^2\right)^2,
$$
where $R=8x_1^2+8x_2^2-x_3^2-x_4^2$.

\section{The unique example in \texorpdfstring{$\Hi^5$}{H5}}\label{h5} 
Here we show how the Gauss parametrization can be used to obtain a
simpler and more direct proof of Tsukada's theorem, which states that
there is a unique rank two curvature homogeneous hypersurface in $\Hi^5$.
We will also recover its basic properties, showing in particular that it
is closely related to both $f_c$'s in the Introduction.

\bigskip

The polar map of such a hypersurface is a surface in the De-Sitter space,
$g:V^2\to \DS^5\subset \R^{5,1}$, i.e.,
$$
\la f,f\ra=-1,\ \ \la f,g\ra=0, \ \ \la g,g\ra=1,\ \ \la \d f,g\ra=0.
$$
By \pref{gp} we know that $\det A_w\neq 0$ is constant for every $w$
in an open subset of~$\Lambda$.

Choose a orthonormal normal frame $\{\xi_0,\xi_1,\xi_2\}$ of
$T^\perp_g V$ with
$-\la \xi_0,\xi_0\ra=\la \xi_1,\xi_1\ra=\la \xi_2,\xi_2\ra=1$.
We call the respective shape operators $A,B,C$ for short, and we can
assume that $\tr B=0$. Write
$w=\cosh(r)\xi_0+\sinh(r)(\cos(t)\xi_1+\sin(t)\xi_2)\in \Lambda$ for
certain $(r,t)\in W\subset\R^2$, $W$ open, and thus $a^2=-\det A_w$ is
constant. Therefore,
$$
a\, A_w=\cosh(r)A+\sinh(r)B_t,\ \ \  B_t=\cos(t)B+\sin(t)C.
$$
Since $a\neq 0$ it easily follows that $A$ is invertible, and hence
\begin{equation}\label{cor6}
\cosh(r)^2\det A + \sinh(r)^2 \det B_t + \cosh(r)\sinh(r)
\tr (A^{-1}B_t)\det A =-1.
\end{equation}
This is equivalent to $\det A=-\det B_t=1$, and $\tr (A^{-1}B_t)=0$.
By \lref{al}, the pair $\{B,C\}$ has the special normal form
\eqref{al1} for $a=c=1$. Since $\tr (A^{-1}B_t)=0$, in this tangent frame
$A$ must have the form $Ae_1=h e_1$, $Ae_2=h^{-1} e_2$, up to a possible
change of the sign of $\xi_0$. Finally, replacing
$\xi_0$ by $\frac{1+h^2}{2h}\xi_0+\frac{1-h^2}{2h}\xi_2$ and
$\xi_2$ by $\frac{1-h^2}{2h}\xi_0+\frac{1+h^2}{2h}\xi_2$,
we can assume that $h=1$ and $A=I$.
We conclude that $V^2$ has constant curvature $1-3a^2$, and that,
in a fixed orthonormal basis $\{e_1,e_2\}$ of $V^2$, the second
fundamental form of $g$ is unique and satisfies
$A_{\xi_0}=a I$, with $A_{\xi_1}$, $A_{\xi_2}$ as in \eqref{vl}
in some orthonormal normal frame that we still call
$\{\xi_0,\xi_1,\xi_2\}$.

We can now easily compute the normal connection 1-forms $w_j^i$, that is,
$\nabla_X^\perp \xi_j=\sum_{i=1}^3w_j^i(X)\xi_i$.
Noticing that
$(-1)^{\delta_0^j}w_j^i+w_i^j=0$, set $w_{i+j}=w_i^j$ for $i<j$.
The Codazzi equations are as usual
$[DA_{\xi_j}]^*=-(-1)^{\delta_0^j}\sum_iw_j^i \circ JA_{\xi_i}$,
where $[DA]=\nabla_{e_1}A(e_2)-\nabla_{e_2}A(e_1)-A[e_1,e_2]$
and $J$ is given by $Je_1=e_2,Je_2=-e_1$. In our case, if
$\beta=\la\nabla_\bullet e_1,e_2\ra$,
$$
w_1\circ JB + w_2\circ JC=0,\ \ \ 
2\beta\circ B = w_1\circ J - w_3\circ JC,\ \ \ 
2\beta\circ C = w_2\circ J + w_3\circ JB,
$$
which determines the normal connection and is independent of $a$.
Using that $-BC=CB=J$ and $B^2=C^2=I$ we easily see that
$w_1=w_2=0$, $w_3=2\beta$. In particular, $\xi_0$ is normal parallel and
$\d w_3=2(3a^2-1)dvol$ since $V^2$ has constant curvature $1-3a^2$.
Furthermore, the Ricci equation implies that
$$
2a^2 = -\la [A_{\xi_1},A_{\xi_2}]e_1,e_2\ra
=-\la R^\perp(e_1,e_2)\xi_1,\xi_2\ra=-\d w_3(e_1,e_2)=2(1-3a^2).
$$
We conclude that $a^2=1/4$, $V^2$ is locally isometric to
$\Sp^2_{1/4}$, and $g$ is unique.

Now, since all shape operators of the minimal Veronese embedding
$g_1:\RP^2_{1/3}\to\Sp^4$ are conjugate to $3^{-1/2}B$, and the shape
operator in the normal parallel direction $\xi_0$ is $I/2$, it is
easy to get an explicit expression for $g$,
$$
g = \frac{1}{\sqrt{3}}(2g_1,1):\RP^2_{1/4}\to\DS^5\subset\R^{5,1}.
$$

Once we computed $g$ we can finally recover $f$. The normal bundle of
$g$ in $\R^{5,1}$ is $\spa\{g, \xi_0=(g_1,2)/\sqrt{3}\,\}\oplus \nu$,
where $\nu$ stands for the normal bundle of $g_1$ and $\nu_1$ its unit
normal bundle. So
$\Lambda=\{c\xi_0+s(\xi,0): \xi\in \nu_1, c^2-s^2=1\}$,
and therefore $f:M^4=\Lambda=\nu_1\times\R\to\Hi^5\subset\R^{5,1}$
is
$$
f(\xi_x,s)=\frac{1}{\sqrt{3}}\left(\cosh(s)g_1(x)
+\sqrt{3}\sinh(s)\xi_x,2\cosh(s)\right).
$$
Since $\nu_1$ as a hypersurface in $\Sp^4$ is
$\SO(3)$-equivariant, so is $f$, with $\SO(3)$ acting on
$\R^5\times\{0\}\subset\R^{5,1}$. Notice that $f$ is clearly complete
since $\nu_1$ is compact and $\|f_*\partial_s\|=1$.
Compare the above expression for $f$ with the one for $f_{-1}$ in
\eqref{f1}.


\end{document}